\documentclass[reqno,10pt]{amsart}
\usepackage{amsmath,amssymb,amsthm,graphicx,mathrsfs,bbm,url}
\usepackage{amsthm}
\usepackage{wrapfig}
\usepackage{enumitem}
\usepackage{mathtools}
\usepackage[all]{xy}
\usepackage[utf8]{inputenc}
\usepackage[usenames,dvipsnames]{color}
\usepackage[colorlinks=true,linkcolor=Red,citecolor=Green]{hyperref}
\usepackage[super]{nth}
\usepackage[open, openlevel=2, depth=3, atend]{bookmark}
\hypersetup{pdfstartview=XYZ}
\usepackage[font=footnotesize]{caption}
\usepackage{fourier}
\usepackage{subcaption}
\usepackage{caption}
\usepackage{tikz}
\usepackage{setspace}
\usepackage[normalem]{ ulem }
\usepackage{soul}
\usetikzlibrary{decorations.pathreplacing}
\usetikzlibrary{arrows}
\usetikzlibrary{shapes.misc}
\usetikzlibrary{shapes.symbols}
\usetikzlibrary{patterns}
\usepackage{epstopdf}
 
\usepackage{hyperref}

\usepackage{listings}
\usepackage{color}

\definecolor{dkgreen}{rgb}{0,0.6,0}
\definecolor{gray}{rgb}{0.5,0.5,0.5}
\definecolor{mauve}{rgb}{0.58,0,0.82}

\theoremstyle{plain}
\newtheorem{theorem}{Theorem}[section]
\newtheorem*{theorem*}{Theorem}

\newtheorem{lemma}[theorem]{Lemma}
\newtheorem{proposition}[theorem]{Proposition}
\newtheorem{corollary}[theorem]{Corollary}

\newtheorem{conjecture}[theorem]{Conjecture}

\theoremstyle{definition}

\theoremstyle{remark}
\newtheorem{remark}[theorem]{Remark}

\numberwithin{equation}{section}

\newcommand{\C}{\mathbb{C}}
\newcommand{\R}{\mathbb{R}}

\newcommand{\Z}{\mathbb{Z}}

\newcommand{\M}{\mathcal{M}}
\newcommand{\N}{\mathbb{N}}

\newcommand{\V}{\mathbb{V}}
\newcommand{\X}{\mathbf{X}}

\newcommand{\HH}{\mathbb{H}}
\newcommand{\Ss}{\mathbb{S}}

\newcommand{\mc}{\mathcal}

\newcommand{\Sym}{\mathrm{Sym}}

\DeclareMathOperator{\E}{\mathcal{E}}

\DeclareMathOperator{\e}{\mathbf{e}}

\DeclareMathOperator{\Id}{\mathrm{Id}}

\newcommand{\be}{\begin{equation}}
\newcommand{\ee}{\end{equation}}

\def\beq{\begin{equation}}
\def\eeq{\end{equation}}

\def \RM{\mathbb{R}}

\def\beq{\begin{equation}}
\def\eeq{\end{equation}}
\def\bea{\begin{eqnarray*}}
\def\eea{\end{eqnarray*}}

% ...

% ...

\title
[Brin's conjecture on frame flow ergodicity]
{Towards Brin's conjecture on frame flow ergodicity: \\ new progress and perspectives}

\author{Mihajlo Ceki\'c}
\address{Institut f\"ur Mathematik, Universit\"at Z\"urich, Winterthurerstrasse 190, CH-8057 Z\"urich, Switzerland}
\email{mihajlo.cekic@math.uzh.ch}

\author{Thibault Lefeuvre}
\address{Université de Paris and Sorbonne Université, CNRS, IMJ-PRG, F-75006 Paris, France.}
\email{tlefeuvre@imj-prg.fr}

\author{Andrei Moroianu}
\address{Université Paris-Saclay, CNRS,  Laboratoire de mathématiques d'Orsay, 91405, Orsay, France}
\email{andrei.moroianu@math.cnrs.fr}

\author{Uwe Semmelmann}
\address{Institut f\"ur Geometrie und Topologie, Fachbereich Mathematik, Universit{\"a}t Stuttgart, Pfaffenwaldring 57, 70569 Stuttgart, Germany
}
\email{uwe.semmelmann@mathematik.uni-stuttgart.de}

\begin{document}

\begin{abstract} 
We report on some recent progress achieved in \cite{Cekic-Lefeuvre-Moroianu-Semmelmann-21} on the ergodicity of the frame flow of negatively-curved Riemannian manifolds. We explain the new ideas leading to ergodicity for nearly $0.25$-pinched manifolds and give perspectives for future work.
\end{abstract}

\maketitle

\section{Ergodicity in dynamical systems}

\subsection{Historical background}

The theory of dynamical systems originates from the study of classical mechanics and the description of solutions to differential equations governing the evolution of points in phase space. A well-known historical example pioneered by Kepler and Newton in the XVII century and later enriched with a modern mathematical language by Poincaré in the late XIX century is our solar system, where planets are identified with points and their motion is governed by the law of gravitation. These mechanical systems, in the absence of a dissipative correction, share the property that they preserve a natural volume form on the phase space.

While it is usually impossible to predict the evolution of a single trajectory
due to an inherent sensitivity to initial conditions, it is tempting to adopt a statistical approach and describe the long-time behaviour of \emph{almost all} points, measured with respect to this natural flow-invariant form. This is the slant of ergodic theory. From this perspective, a natural property one may investigate on a given dynamical system is the \emph{equidistribution} of a generic point in phase space, namely, whether it will spend in each region of the phase space an average time proportional to its volume. Phrased in mathematical language, ergodicity is the property that any measurable subset that is invariant by the dynamical transformation must have zero or full measure.

These physical considerations paved the way for a more systematic search of ergodic dynamical systems in mathematics. In a seminal article \cite{Hopf-36}, using what is now known as the classical \emph{Hopf argument}, Hopf proved that geodesic flows on closed negatively-curved surfaces are ergodic with respect to a natural smooth measure called the Liouville measure, providing one of the first rigorous examples of chaotic systems of geometric flavour. Later, Anosov \cite{Anosov-67} introduced in his thesis the notion of \emph{uniformly hyperbolic flows} (also known as \emph{Anosov flows} nowadays) and showed that they are ergodic whenever they preserve a smooth measure. Moreover, he proved that all geodesic flows on negatively-curved Riemannian manifolds are uniformly hyperbolic and thus ergodic. From a statistical perspective, these flows are now well understood and finer properties such as \emph{mixing} or even \emph{exponential mixing} are (almost) completely settled, see Liverani \cite{Liverani-04} and Tsujii-Zhang \cite{Tsuji-Zhang-20}, among other references on this question.

Shortly after Anosov's work, Brin-Pesin \cite{Brin-Pesin-74}, Pugh-Shub \cite{Pugh-Shub-00}, and others, investigated more general systems exhibiting a weaker form of hyperbolic behaviour, known as \emph{partially hyperbolic systems}. While these dynamics still preserve some expanding and contracting directions, they also come with other \emph{neutral} or \emph{central} directions, in which the map/flow may behave infinitesimally as an isometry for instance. Historical examples of partially hyperbolic dynamics are provided by \emph{frame flows} over closed negatively-curved Riemannian manifolds, which are at the core of the present article.

\subsection{Ergodicity of the frame flow. Statement of results}

\label{ssection:intro2}

Let $(M^n,g)$ be a closed connected oriented $n$-dimensional Riemannian manifold with negative sectional curvature. Let $SM := \left\{(x,v) \in TM ~|~ |v|_g = 1\right\}$ be the unit tangent bundle over $M$ and let 
\[
FM := \big\{(x,v, \e_2, \dotso, \e_n) ~|~   (x, v) \in SM,\,    (v, \e_2, \dotso, \e_n)  \, \mathrm{oriented\,\,orthonormal\,\,basis\,\,of\,\,} T_xM\big\}
\]
be the principal $\mathrm{SO}(n)$-bundle (over $M$) of oriented orthonormal frames.

We can also consider $p : FM \to SM, (x, v, \e_2, \dotso, \e_n) \mapsto (x, v)$ as a principal $\mathrm{SO}(n-1)$-bundle over $SM$, that is, a point $w \in FM$ over $(x,v) = p(w)$ corresponds to an orthonormal frame $(\e_2,\dotso ,\e_n)$ of the orthogonal complement $v^\perp \subset T_xM$. Denoting by $\nabla$ the Levi-Civita connection on $M$, the geodesic flow $(\varphi_t)_{t \in \R}$ is defined on $SM$ by setting for $t\in \R, (x,v) \in SM$, $\varphi_t(x,v) := (\gamma_{x,v}(t), \dot{\gamma}_{x,v}(t))$, where $t \mapsto \gamma_{x,v}(t)$ is the unit-speed curve on $M$ solving the geodesic equation $\nabla_{\dot{\gamma}_{x,v}}\dot{\gamma}_{x,v} = 0$ with initial conditions $\gamma_{x,v}(0)=x, \dot{\gamma}_{x,v}(0)=v$.

The \emph{frame flow} $(\Phi_t)_{t \in \R}$ on $FM$ is then defined as follows: 
\[
\Phi_t(x,v,\e_2,...,\e_n) := (\varphi_t(x,v), \tau_{x,v}(t)\e_2, \dotso, \tau_{x,v}(t)\e_n),\]
where $\tau_{x,v}(t)$ denotes the parallel transport along the geodesic segment $\gamma_{x,v}([0,t])$ with respect to the connection $\nabla$. Moreover, since the geodesic flow $(\varphi_t)_{t \in \R}$ preserves the \emph{Liouville measure} $\mu$ on $SM$, the frame flow $(\Phi_t)_{t \in \R}$ preserves a smooth measure $\omega$ on $P$ induced by $\mu$ and the Haar measure on $\mathrm{SO}(n-1)$. It is therefore natural to study the ergodicity of the frame flow with respect to the smooth measure $\omega$.

It was first shown by Brin \cite{Brin-75-1} (for $n=3$) and later by Brin-Gromov \cite{Brin-Gromov-80} (for $n$ odd and different from $7$) that negatively-curved $n$-dimensional manifolds have an ergodic frame flow. As we shall see below in \S\ref{section:topology}, once the dynamical framework is settled, the proof boils down to a (non-trivial) statement in algebraic topology on the classification of topological structures over even dimensional spheres. It is however hopeless to expect \emph{all} negatively-curved manifolds to have an ergodic frame flow: indeed, it can be checked that Kähler manifolds of real dimension $n=2m \geq 4$ such as compact quotients $\Gamma \backslash \C\HH^m$ of the complex hyperbolic space (where $\Gamma \leqslant \mathrm{Isom}(\C\HH^m)$ is a lattice) do not have an ergodic frame flow\footnote{This may be seen as follows: the complex structure $J$ of a K\"ahler manifold commutes with parallel transport $\tau_{x, v}(t)$, so the set $\left\{(x, v, \e_2, \dotso, \e_n) \in FM \mid g_x(Jv, \e_2) \geq 0\right\}$ is invariant and has positive, but not full measure. In the even-dimensional case, the situation therefore requires additional care.} due to the reduction of their holonomy group from $\mathrm{SO}(n)$ to $\mathrm{U}(m)$. 

Denoting by $\kappa_g(u \wedge v)$ the sectional curvature of the $2$-plane spanned by $u,v \in TM$, we will say that $(M,g)$ has \emph{$\delta$-pinched negative curvature} for some $\delta \in (0,1]$ if there exists a constant $C > 0$ such that the following inequalities hold:
\[
-C \leq \kappa_g(u \wedge v) \leq -C \delta.
\]
Note that, up to rescaling the metric, one can always assume that $C=1$. Since Kähler manifolds are at most $0.25$-pinched by a result of Berger \cite{Berger-60-1}, Brin \cite{Brin-82} stated the following natural conjecture:

\begin{conjecture}[Brin '82]
\label{conjecture:brin}
If $(M,g)$ is $\delta$-pinched for some $\delta > 0.25$, then the frame flow is ergodic.
\end{conjecture}

More generally, Brin conjectures in the same article (see \cite[Conjecture 2.9]{Brin-82}) that the frame flow should be ergodic as long as there is no reduction of the holonomy group of the manifold. However, up to now, ergodicity of the frame flow in dimension $7$ and on even-dimensional manifolds was only known for nearly-hyperbolic manifolds, that is, manifolds with a pinching $\delta$ very close to $1$: strictly greater than $0.8649...$ in even dimensions different from $8$, due to Brin-Karcher \cite{Brin-Karcher-83}, and strictly greater than $0.9805...$ in dimensions $7$ and $8$, due to Burns-Pollicott \cite{Burns-Pollicott-03}. There has been no progress on Conjecture \ref{conjecture:brin} in the past twenty years, until our result in \cite{Cekic-Lefeuvre-Moroianu-Semmelmann-21}:

\begin{theorem}[Ceki\'c-Lefeuvre-Moroianu-Semmelmann '21]
\label{theorem:main1}
The frame flow is ergodic if the manifold is $\delta(n)$-pinched, where $\delta(7) \sim 0.497$, and asymptotically $\delta(n) \sim 0.277$ for $n \equiv 2 \mod 4$ and $\delta(n) \sim 0.557$ for $n \equiv 0 \mod 4$.
\end{theorem}

\begin{center}
\begin{figure}[htbp!]
\includegraphics[scale=0.47]{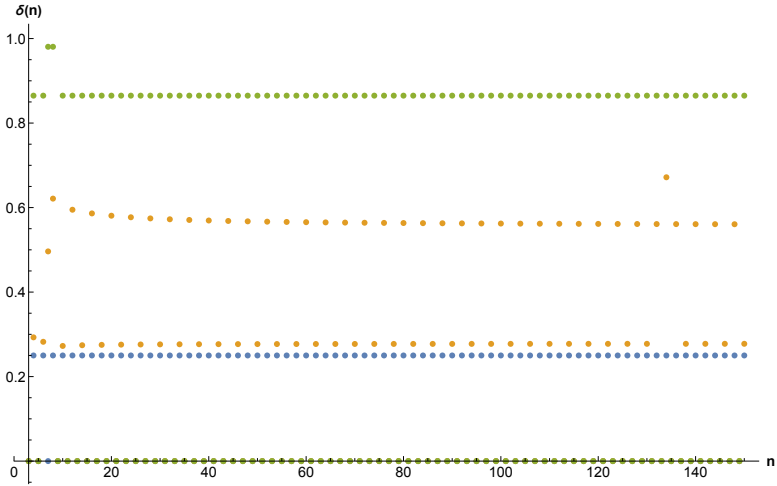}
\caption{In green: the bounds existing in the literature \cite{Brin-Gromov-80,Brin-Karcher-83,Burns-Pollicott-03}. In orange: the bounds provided by Theorem \ref{theorem:main1}. In blue: the conjectural $0.25$ threshold.}
\label{figure}
\end{figure}
\end{center}

The precise version of Theorem \ref{theorem:main1} can be found in \cite[Theorem 1.2]{Cekic-Lefeuvre-Moroianu-Semmelmann-21}. Nevertheless, the numerical value of $\delta(n)$ is depicted in Figure \ref{figure} for $n \in \left\{4,...,150\right\}$. 

While former attempts to prove Conjecture \ref{conjecture:brin} were mostly based on algebraic topology or on the geometry of the universal cover of the manifold, the strategy of proof for Theorem \ref{theorem:main1} is different and relies on the introduction of new ideas from Riemannian geometry. More precisely, we make systematic use of the twisted version of an energy identity on the unit tangent bundle known as the \emph{Pestov identity}, first introduced by Mukhometov \cite{Mukhometov-75,Mukhometov-81} and Amirov \cite{Amirov-86}, then in its classical form by Pestov and Sharafutdinov \cite{Pestov-Sharafutdinov-88,Sharafutdinov-94} and finally stated in full generality by Guillarmou-Paternain-Salo-Uhlmann \cite{Guillarmou-Paternain-Salo-Uhlmann-16}. This identity has found several applications in the past twenty years, be it in the study of inverse spectral problems, see Croke-Sharafutdinov \cite{Croke-Sharafutdinov-98}, or in tensor tomography \cite{Paternain-Salo-Uhlmann-13}.

Our argument actually consists of three distinct parts, each of them belonging to a different area of mathematics:

 \begin{enumerate}
 \item \textbf{Hyperbolic dynamical systems:} the frame bundle $p: FM \to SM$ is a principal $\mathrm{SO}(n-1)$-bundle; using the hyperbolic structure of the geodesic flow, one can prove that the non-ergodicity of the frame flow entails the existence of a subgroup $H \lneqq \mathrm{SO}(n-1)$ called the \emph{transitivity group} and an $H$-principal subbundle $FM \supset Q \to SM$ which is an ergodic component of the frame flow, see \S\ref{section:dynamics}.
 \item \textbf{Algebraic topology:} in particular, restricting to an arbitrary point $x_0 \in M$, one obtains an $H$-subbundle $Q_{x_0} \to S_{x_0}M \simeq \Ss^{n-1}$ of the frame bundle of $\mathbb{S}^{n - 1}$; this is already a strong topological constraint, called \emph{reduction of the structure group} of the sphere, which rules out most possible cases for the subgroup $H$, see \S\ref{section:topology}.
 \item \textbf{Riemannian geometry:} when topology is not sufficient to rule out the existence of a structure group reduction, we show that, according to the possible values of $H$, one can produce a smooth flow-invariant section $f \in C^\infty(SM, \E)$, where $\E = \pi^* \Lambda^p TM$ (with $p=1,2,3$) or $\E=\pi^* \mathrm{Sym}^2 TM$, and $\pi : SM \to M$ is the projection. In turn, the existence of such object can be ruled out by means of the \emph{twisted Pestov identity} whenever the pinching $\delta$ is sufficiently large, see \S\ref{ssection:pestov}.
 \end{enumerate}

Three exotic dimensions appear in this setting: $n=7,8,$ and $134$. They correspond to special topological structures possibly carried by the spheres $\Ss^6, \Ss^7,$ and $\Ss^{133}$, respectively. The induced flow-invariant sections obtained in point (3) above then take values in $\pi^*\Lambda^pTM$ for $p=2,3$. It is more difficult to rule out the existence of such objects and our method requires a larger pinching than in other cases, see Figure \ref{figure} (for instance, the orange dot on the right-hand side corresponds to the case $n=134$).

On the other hand, when $n \equiv 2 \mod 4$, the maximal number of linearly independent vector fields on the sphere $\mathbb{S}^{n - 1}$ is 1, which simplifies our analysis and eventually yields a flow-invariant section of $\pi^*TM$, whereas in the case $n \equiv 0 \mod 4$, it is at least 3 and we get a flow-invariant orthogonal projector which is a section of $\pi^*\Sym^2 TM$.
 
The purpose of this article is to explain the circle of ideas leading to Theorem \ref{theorem:main1} and unlocking the long-standing problem of frame flow ergodicity.
This strategy may also prove seminal for additional progress. It comes with many interesting new open questions, summed up in the last paragraph \S\ref{section:questions}, some of them being related to (polynomial) structures over spheres. Eventually, it would be interesting to understand to what extent the techniques of the present paper may apply to the broader setting of non-isometric partially hyperbolic dynamics. This is left for future investigation. \\

\noindent \textbf{Acknowledgment:} M.C. has received funding from an Ambizione grant (project number
201806) from the Swiss National Science Foundation. We are grateful to an anonymous referee for giving us a historical overview of the Pestov identity.

\section{Principal bundle extensions of Anosov flows}

\label{section:dynamics}

\subsection{A partially hyperbolic flow}

Let $\M$ be a smooth closed manifold. We recall that a vector field $X \in C^\infty(\M,T\M)$ generates an \emph{Anosov flow} $(\varphi_t)_{t \in \R}$ if there exists a continuous flow-invariant splitting of the tangent bundle $T\M = \R X \oplus E^s_\M \oplus E^u_\M$ into flow-direction, stable and unstable bundles, and uniform constants $C,\lambda > 0$ such that for all $t \geq 0$:
\begin{equation}
\label{equation:anosov}
 \|d\varphi_{t} v\| \leq C e^{-\lambda t} \|v\|,\ \forall v \in E^{s}_\M,  \qquad  \|d\varphi_{- t} v\| \leq C e^{-\lambda t} \|v\|,\ \forall v \in E^{u}_\M,
\end{equation}
where $\|\bullet\|$ is the norm induced by an arbitrary Riemannian metric on $\mc{M}$. In the following, we will further assume that $X$ preserves a \emph{smooth measure} $\mu$. In particular, this implies that $(\varphi_t)_{t \in \R}$ is ergodic with respect to $\mu$. The goal of this paragraph is to study the ergodic properties of some specific \emph{extensions} of $(\varphi_t)_{t \in \R}$ which we now describe. 

Given $p : P \to \M$, a $G$-principal bundle, where $G$ is a compact Lie group, we say that a flow $(\Phi_t)_{t \in \R}$ is a \emph{principal extension} of $(\varphi_t)_{t \in \R}$ to the bundle $P$ if it satisfies the following two conditions:
\begin{equation}
\varphi_t \circ p = p \circ \Phi_t, \qquad R_g \circ \Phi_t = \Phi_t \circ R_g, \qquad \forall t \in \R,\ \forall g \in G,
\end{equation}
where $R_g : P \to P$ denotes the fiberwise right-action of the group. Such a flow then preserves a natural smooth measure $\omega$ which can be locally written as $\omega = \mu \times \mathrm{Haar}_G$, the right factor being the normalized Haar measure on the group.

Understanding the ergodicity of $(\Phi_t)_{t \in \R}$ with respect to $\omega$ is a very natural question in order to describe the long-time statistical properties of the extended flow. As mentioned in \S\ref{ssection:intro2}, an archetypal example fitting in this framework is the frame flow over a negatively-curved Riemannian manifold $(M, g)$, which will be further discussed in \S\ref{section:topology} (in this case $P=FM$ and $\mc{M} = SM$).

Principal extensions of Anosov flows are intrinsically more complicated to study due to their lack of uniform hyperbolicity. Indeed, the \emph{vertical direction} $\V := \ker d p$ now becomes a \emph{neutral} or \emph{central} direction, in the sense that the differential of the flow $(\Phi_t)_{t \in \R}$ acts as a linear isometry on $\V$ and the tangent bundle to $P$ then splits as
\[
TP = \R X_P \oplus E^s_P \oplus E^u_P \oplus \V,
\]
where $X_P$ is the vector field generating $(\Phi_t)_{t \in \R}$ and $E^{s,u}_P$ satisfy an expanding/contracting property similar to \eqref{equation:anosov}. Note that the subbundles $E^{s,u}_P$ also integrate to produce a (Hölder-continuous) foliation on $P$ by strong stable and unstable manifolds $W^{s,u}_P$, see Pesin \cite{Pesin-04} or Hasselblatt-Pesin \cite{Hasselblatt-Pesin-06} for the related case of partially hyperbolic diffeomorphisms.

\subsection{Transitivity group. Parry's free monoid}

Following Hopf's argument in the Anosov case, it is natural to expect (at least heuristically) that the ergodic component of an arbitrary point $z \in P$ consists of all the other points $z' \in P$ that one can reach from $z$ by following a concatenation of flow- and so-called us-paths, namely, paths that are either fully contained in a flowline of $(\Phi_t)_{t \in \R}$ or in a strong stable/unstable leaf $W^{s,u}_P$. The \emph{full accessibility} of a flow is the property that any other point $z' \in P$ can be reached from $z$ by such a concatenation of paths and it is expected that volume-preserving partially hyperbolic dynamical systems are ergodic whenever they are accessible\footnote{A refinement of this notion is the \emph{essential accessiblity}, that is, accessibility up to measure zero, but this will not be needed here.}: this is known as the Pugh-Shub conjecture \cite{Pugh-Shub-00}. Under a certain additional center bunching assumption, the Pugh-Shub conjecture was proved by Burns-Wilkinson \cite{Burns-Wilkinson-10}. Taking advantage of the very algebraic structure of principal extensions of Anosov flows, Brin \cite{Brin-75-1} translated the accessibility property into a key algebraic notion, called the \emph{transitivity group}: this is a subgroup $H \leqslant G$ (well-defined up to conjugacy in $G$) describing all the points in a fiber that are reachable by flow- and us-paths, which we now describe.

\begin{center}
\begin{figure}[htbp!]
\includegraphics[scale=1]{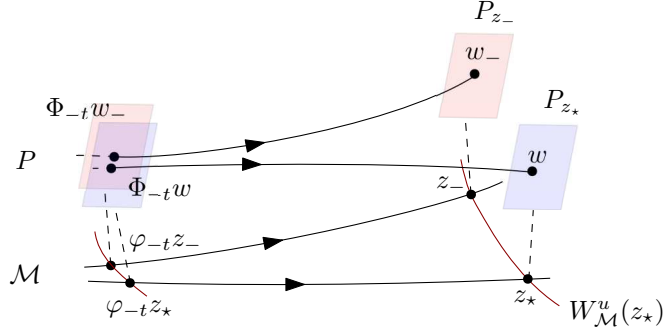}
\caption{Bottom: black lines represent flowlines of $(\varphi_t)_{t \in \R}$ on $\M$; red lines represent strong unstable leaves. Top: black lines represent flowlines of $(\Phi_t)_{t \in \R}$ on $P$; blue and red parallelograms represent fibers of $P$.}
\label{figure2}
\end{figure}
\end{center}

It will be convenient to fix an arbitrary periodic point $z_\star \in \M$ for the flow $(\varphi_t)_{t \in \R}$, generating a periodic orbit $\gamma_\star \subset \M$ of period $T_\star$. Denote by $\mc{H}$ the set of orbits of $(\varphi_t)_{t \in \R}$ that are \emph{homoclinic} to $\gamma_\star$, namely, which accumulate in the past and in the future to $\gamma_\star$. Volume preserving (and more generally, transitive) Anosov flows satisfy that $\mc{H}$ is dense in $\M$. Given $\gamma \in \mc{H}$ and a point $w$ in the fiber $P_\star := P_{z_\star}$ (which we can think of as a frame for instance), there is a natural way to ``parallel transport'' $w$ \emph{along} $\gamma$ (even though $\gamma$ has infinite length!) in order to produce another point, denoted by $\rho(\gamma)w$. This construction goes as follows (see Figure \ref{figure2}, and \cite{Cekic-Lefeuvre-21-1} for more details): 
\begin{enumerate}
\item One picks an arbitrary point $z_- \in \gamma \cap W^u_\M(z_\star)$ (i.e. such that $d_\M(\varphi_{-t}z_\star,\varphi_{-t}z_-) \to 0$ as ${t \to +\infty} $, and this convergence is exponentially fast); then, in the fiber $P_{z_-}$ over $z_-$, there exists a \emph{unique} point $w_- \in W^u_P(w)$ such that $d_P(\Phi_{-t}w, \Phi_{-t}w_-) \to 0$ as $t \to +\infty$. The map $P_\star \to P_{z_-}, w \mapsto w_-$ is called the \emph{unstable holonomy}.
\item One then ``pushes'' $w_-$ by the flow $(\Phi_t)_{t \in \R}$ until it reaches a point $w_+ := \Phi_T(w_-)$ over $z_+ \in \gamma \cap W^s_\M(z_\star)$, where $T > 0$ is the unique time such that $z_+ = \varphi_T(z_-)$;
\item Eventually, applying a similar (but stable this time) holonomy to (1), one can produce an element $\rho(\gamma) w \in P_\star$.
\end{enumerate}

It is important to have in mind that the above-mentioned holonomies are actually parallel transports in strong stable/unstable/flow leaves of $(\varphi_t)_{t \in \R}$ with respect to a \emph{dynamical} $G$-principal connection $\nabla$ on the principal bundle $P \to \M$ induced by the lifted flow $(\Phi_t)_{t \in \R}$. If we identify the fiber $P_\star \simeq G$, then $\rho(\gamma)$ can be identified with left multiplication by an element of the group $G$ itself, that is, $\rho(\gamma) \in G$. Hence, introducing Parry's free monoid\footnote{A monoid is a set endowed with an associative product, a neutral element, but no inverse.} $\mathbf{G}$ as the formal set of words
 \[
 \mathbf{G} := \left\{ \gamma_{i_1}^{k_1} ... \gamma_{i_p}^{k_p} \mid p \in \N, k_j \in \N, \gamma_{i_j} \in \mc{H}, j = 1, ..., p\right\},
 \]
 we see that the above-mentioned construction produces a natural representation $\rho : \mathbf{G} \to G$ whose image $H := \overline{\rho(\mathbf{G})} \leqslant G$ is called the \emph{transitivity group}. Note that this subgroup is only well-defined up to conjugacy as it requires to choose a (non-canonical) identification $P_\star \simeq G$. Moreover, it can be checked that $H$ is nothing but the \emph{holonomy group} of the dynamical connection $\nabla$ on $P$, see \cite{Lefeuvre-21} for instance. This object turns out to be crucial in understanding the ergodicity of the flow $(\Phi_t)_{t \in \R}$:
 
 \begin{proposition}[Brin]\label{prop:Brin1}
 The flow $(\Phi_t)_{t \in \R}$ is ergodic if and only if $H = G$.
 \end{proposition}

Following Proposition \ref{prop:Brin1}, a sound strategy to prove ergodicity is therefore to assume that $H \lneqq G$ is a strict subgroup and to seek a contradiction. The key idea is to show that whenever $H \neq G$, one can produce additional flow-invariant geometric structures over $\M$ using representation theory and to prove by means of topological and/or geometric arguments that actually such structures cannot exist.

\subsection{Flow-invariant structures}

\label{ssection:invariant}


Let us recall at this stage a general principle of differential geometry one has to bear in mind:
\begin{equation}\label{principle:P}
\begin{minipage}{0.85\textwidth}
	\emph{``Any object that is algebraically $H$-invariant will produce a global smooth flow-invariant object over the base manifold $\M$}.''
\end{minipage}\tag{P}
\end{equation}
In addition, this object will be parallel with respect to the dynamical connection. However, the dynamical connection is only H\"older regular, and to make Principle \eqref{principle:P} rigorous one needs non-Abelian Liv\v sic theory, developed in \cite{Cekic-Lefeuvre-21-1}. This idea has two important manifestations. First of all, if $\rho : G \to \mathrm{Aut}(V)$ is a finite-dimensional representation whose restriction to $H$ fixes a non-zero vector $\xi \in V$, then the corresponding associated vector bundle $P \times_\rho V \to \M$ carries a \emph{parallel section} with respect to the connection induced by the dynamical connection $\nabla$, and this object will have a constant algebraic type\footnote{For instance, if $\xi$ is a skew-symmetric endomorphism squaring to $-\mathbbm{1}$, then the parallel section will have the same properties.}. Secondly, $H$ is invariant by the action by itself (as a subgroup of $G$), and \eqref{principle:P} leads to the following remarkable result:
 
 \begin{proposition}[Brin]\label{prop:Brin2}
 There exists a flow-invariant $H$-principal subbundle $Q \subset P$ over $\M$.
 \end{proposition}
 
Non-ergodicity of $(\Phi_t)_{t \in \R}$ thus entails the existence of a \emph{strict} subbundle $Q \to \M$ of the $G$-principal bundle $P \to \M$. This is already a strong topological constraint called \emph{reduction of the structure group of the bundle}.

\section{Structure group reduction of the frame bundle on spheres}

\label{section:topology}

We now specialize the considerations above to the case where $(M^n,g)$ is a smooth compact oriented Riemannian manifold with negative sectional curvature of dimension $n \geq 3$, $\M = SM$ is the unit tangent bundle, $X$ the geodesic vector field, and $P = FM \to SM$ the principal $\mathrm{SO}(n-1)$-bundles of frames as explained in \S\ref{ssection:intro2}, on which we study the frame flow $(\Phi_t)_{t \in \R}$.

%
%
%
Assuming that the frame flow is not ergodic, the transitivity group $H$ is a strict subgroup of $\mathrm{SO}(n-1)$ and by Proposition \ref{prop:Brin2}, there exists an $H$-principal bundle $Q \subset FM$ over $SM$. Restricting $Q$ to the fiber $S_xM\simeq \mathbb{S}^{n-1}$ over some $x\in M$ defines an $H$-principal bundle $Q_x\subset F(\mathbb{S}^{n-1})$ over $S_xM$, i.e. a \emph{structure group reduction} of the orthonormal frame bundle of the round sphere to $H$. Alternatively, this can be seen as follows: since the two hemispheres of $\Ss^{n-1}$ are contractible (hence, any bundle over these is trivial), an $\mathrm{SO}(n-1)$-principal bundle over $\Ss^{n-1}$ is simply given by the data of a \emph{clutching function} $c$ at the equator $c : \Ss^{n-2} \to \mathrm{SO}(n-1)$, taking values in $\mathrm{SO}(n-1)$, and defined up to homotopy (equivalently, $c$ is an element of the homotopy group $\pi_{n-2}(\mathrm{SO}(n-1))$). The bundle admits a reduction of its structure group to $H \leqslant \mathrm{SO}(n-1)$ if $c$ can be factored through $H$ (here $\iota : H \to \mathrm{SO}(n-1)$ is the embedding):
\[
\xymatrix{
     & H \ar[d]^{\iota} \\
    \Ss^{n-2} \ar[ur]^{c_0} \ar[r]_c & \mathrm{SO}(n-1)
  }
\]
This fact alone has already strong topological consequences. Indeed, using the work of Adams \cite{Adams-62}, Leonard \cite{Leonard-71}, and \v Cadek--Crabb \cite{Cadek-Crabb-06}, one can prove:

\begin{proposition}\label{prop:red} 
The following holds:
\begin{enumerate}
\item If $n\ge 3$ is odd and $n\ne 7$, there is no reduction of the structure group of $\mathbb{S}^{n-1}$ to a strict subgroup of $\mathrm{SO}(n-1)$.

\item For all other $n\ge 4$, if there exists a reduction of the structure group of $\mathbb{S}^{n-1}$ to a strict subgroup $H$ of $\mathrm{SO}(n-1)$, then up to conjugation, $H$ is contained in one of the following subgroups $K$ of $\mathrm{SO}(n-1)$:

\begin{itemize}
\item If $n=7$, $K= \mathrm{U}(3)\subset \mathrm{SO}(6)$;
\item If $n=8$, $K = \mathrm{G}_2$ or $K=\mathrm{SO}(p) \times \mathrm{SO}(7-p)\subset \mathrm{SO}(7)$ with $p = 1,2,3$;
\item If $n=134$, $K=\mathrm{E}_7\subset \mathrm{SO}(133)$ or $K = \mathrm{SO}(132)\subset \mathrm{SO}(133)$;
\item If $n \equiv 2 \mod 4$, $n \neq 134$, $K = \mathrm{SO}(n-2)\subset \mathrm{SO}(n-1)$;
\item If $n \equiv 0 \mod 4$, $n \neq 8$, $K = \mathrm{SO}(p) \times \mathrm{SO}(n-1-p)\subset \mathrm{SO}(n-1)$ with $1 \leq p \leq (n-2)/2$.
\end{itemize}
\end{enumerate}
\end{proposition}

Note that the Brin--Gromov result \cite{Brin-Gromov-80} about the ergodicity of the frame flow on negatively curved compact Riemannian manifolds in odd dimensions different from $7$ is a direct consequence of Proposition \ref{prop:Brin1} and Proposition \ref{prop:red} (1) (which was proved by Leonard \cite{Leonard-71}). 

\begin{remark}
\label{remark}
Let us also point out that, unlike other topological reductions which \emph{do} appear, we actually do not know whether the $\mathrm{E}_7$-structure on $\Ss^{133}$ exists or not. This is still an open question.
\end{remark}

We will now use the discussion of \S\ref{ssection:invariant} in order to produce new flow-invariant geometric objects whenever the flow is not ergodic, that is, whenever the transitivity group $H$ is strictly contained in $\mathrm{SO}(n-1)$. There is a natural associated vector bundle $\mc{V} = FM \times_{\rho} \R^{n-1} \to SM$, given by the canonical representation $\rho : \mathrm{SO}(n-1) \to \mathrm{Aut}(\R^{n-1})$, called the \emph{normal bundle}. This bundle is also isomorphic to the \emph{vertical} bundle of the spherical fibration $SM \to M$, that is, the vector bundle whose fiber at $(x,v) \in SM$ is the $(n-1)$-dimensional space $v^\perp \subset T_xM$. There is a natural way to parallel transport sections of this bundle along geodesic flowlines with respect to the (lift of the) Levi-Civita connection and, therefore, it makes sense to talk about flow-invariant sections.

The key point is then that for each group $K$ occurring in Proposition \ref{prop:red} (2), one can find non-zero $K$-invariant vectors in some tensorial representations. More precisely:
\begin{itemize}
\item $\mathrm{U}(3)\subset \mathrm{SO}(6)$ preserves a non-zero 2-form in $\Lambda^2\mathbb{R}^6$;
\item $\mathrm{G}_2$ preserves a non-zero $3$-form in $\Lambda^3\mathbb{R}^7$;
\item $\mathrm{E}_7\subset \mathrm{SO}(133)$ preserves a non-zero $3$-form\footnote{Indeed, the embedding of $\mathrm{E}_7$ in $\mathrm{SO}(133)$ is obtained via the adjoint representation of $\mathrm{E}_7$ on its Lie algebra $\mathfrak{e}_7=\mathbb{R}^{133}$, so $\mathrm{E}_7$ preserves the canonical 3-form of $\mathfrak{e}_7$.} in $\Lambda^3\mathbb{R}^{133}$;
\item $\mathrm{SO}(n-2)\subset \mathrm{SO}(n-1)$ preserves a unit vector in $\mathbb{R}^{n-1}$;
\item For $1\le p\le (n-2)/2$, $\mathrm{SO}(p) \times \mathrm{SO}(n-1-p)\subset \mathrm{SO}(n-1)$ preserves the orthogonal projection of $\mathbb{R}^{n-1}$ onto $\mathbb{R}^{p}$.
\end{itemize}
Following \S\ref{ssection:invariant}, this implies in turn that in all these cases, one can produce a smooth parallel (with respect to the dynamical connection) object over $SM$. In particular, this object is flow-invariant:

\begin{theorem}\label{fis}
If the frame flow of $M$ is not ergodic, there exists a non-vanishing flow-invariant section $f \in C^\infty(SM,\E)$, where $\E$ is one of the following bundles: 
\[
\begin{array}{ll}
  \text{(1) } \E = \mc{V}, & \qquad \text{(2) } \E = \Lambda^2 \mc{V} \text{ for $n=7$}, \\
 \text{(3) } \E = \Lambda^3 \mc{V} \text{ for $n=8$ or $n=134$}, &   \qquad \text{(4) } \E = \mathrm{Sym}^2 \mc{V}.
\end{array}
\]
\end{theorem}

For instance, if $n=7$, $H \leqslant \mathrm{U}(3)$ and thus $H$ fixes a an element $\Lambda^2 \R^6$ which is an almost complex structure, namely, a skew-symmetric endomorphism squaring to $-\mathbbm{1}_{\R^6}$. In turn, this gives rise to the existence of a flow-invariant section $f \in C^\infty(SM,\Lambda^2\mc{V})$ of constant algebraic type, that is, such that for every $(x,v) \in SM$, $f(x,v)$ is a skew-symmetric endomorphism on $v^\perp$ squaring to $-\mathbbm{1}_{v^\perp}$. Our aim is now to show that, under some pinching condition, one can rule out the existence of such a flow-invariant geometric structure on the unit tangent bundle $SM$.

\section{Bounding the degree of flow-invariant objects}

\label{section:geometry}

The existence of flow-invariant structures provided by Theorem \ref{fis} is not enough to obtain a contradiction. Indeed, such objects do actually exist in some settings as we shall see below in \S\ref{ssection:examples}. However, under some pinching condition and using Fourier analysis on $SM$, one can actually describe very accurately the analytic properties of flow-invariant objects. This is eventually what will give us the contradiction we are seeking.

\subsection{Fourier degree of sections}

If $\Delta$ denotes the Laplacian acting on functions on the round sphere $\mathbb{S}^{n-1}\subset \mathbb{R}^n$, the eigenspaces
$$
\Omega_k:=\left\{f\in C^\infty(\mathbb{S}^{n-1}) \mid \Delta f=k(n+k-2)f\right\}, \quad k \geq 0,
$$
consist in restrictions to $\mathbb{S}^{n-1}$ of harmonic homogeneous polynomials of degree $k$ on $\mathbb{R}^n$. Elements in $\Omega_k$ are called spherical harmonics of degree $k$. 
Every function $f\in C^\infty(\mathbb{S}^{n-1})$ has a unique decomposition (in the $L^2$-sense) as $f=\sum_{j\ge 0}f_j$, with $f_j\in\Omega_j$. This decomposition also applies to functions defined on the sphere bundle $SM$ of a Riemannian manifold $(M,g)$, or more generally to sections of the pull-back to $SM$ of vector bundles over $M$. 

More precisely, if $\E:=\pi^*{E}$ denotes the pull-back to $SM$ of a vector bundle $E$ over $M$, its restriction to any fiber $S_xM\simeq \mathbb{S}^{n-1}$ is trivial, so the restriction of any section $f\in C^\infty(SM,\E)$ to $S_xM$ can be identified with a vector-valued function $f|_{S_xM}:\mathbb{S}^{n-1}\to {E}_x$. The vertical Laplacian $\Delta_{\mc{V}}$ acts on sections of $\E$ and satisfies $(\Delta_{\mc{V}}f)|_{S_xM}=\Delta(f|_{S_xM})$ for every $x\in M$. Correspondingly, setting for $x \in M$,
$$
\Omega_k(\E)_x :=\left\{f\in C^\infty(S_xM,\E) \mid \Delta_{\mc{V}} f =k(n+k-2)f\right\},
$$
we get a vector bundle $\Omega_k(\E) \to M$ and the decomposition of any section $f\in C^\infty(SM,\E)$ as $f=\sum_{j\ge 0}f_j$, with $f_j\in C^\infty(M,\Omega_j(\E))$. If the above sum is finite, i.e. $f=\sum_{j=0}^{k}f_j$ with $f_k\ne 0$, we say that $f$ has \emph{finite degree} $k$. If the above sum only contains even (resp. odd) spherical harmonics, i.e. $f = \sum_{j \geq 0} f_{2j}$ (resp.  $f = \sum_{j \geq 0} f_{2j+1}$), we say that $f$ is \emph{even} (resp. \emph{odd}).

From the description of $\Omega_k$ as the set of harmonic homogeneous polynomials on $\RM^n$ of degree $k$, it easily follows that $\Omega_k(\E)_x$ can be identified with $\Sym^k_0(T^*_xM)\otimes {E}_x$ by the tautological map $\pi^*: \Sym^k_0(T^*_xM)\otimes {E}_x \to \Omega_k(\E)_x$ defined as follows: if $K\in  \Sym^k_0(T^*_xM)$ is a trace-free symmetric tensor of degree $k$ and $s\in {E}_x$, one defines 
\begin{equation}
\label{identif}\pi^*(K\otimes s)_{(x,v)}:=\frac1{k!}K(v,\ldots,v)s_x,\qquad \forall v\in S_xM.
\end{equation}
More generally, \eqref{identif} identifies $\Sym^k(T^*M) \otimes {E}$ with $\oplus_{j \geq 0} \Omega_{k-2j}(\E)$ with the convention that $\Omega_j(\E) = \left\{0\right\}$ for $j < 0$.

Whenever $E$ is equipped with a metric connection $\nabla^{E}$, we can consider the pull-back connection $\pi^*\nabla^{E}$ on $\E:=\pi^*E$. We set $\X := (\pi^*\nabla^{E})_X$, where $X$ is the geodesic vector field on $SM$, which is nothing but the generator of the parallel transport of sections of $E$ with respect to $\nabla^{E}$ along geodesics. This operator may be seen to have the mapping property
\begin{equation}
\label{eq:XX}
\X : C^\infty(M,\Omega_k(\E)) \to C^\infty(M,\Omega_{k-1}(\E)) \oplus C^\infty(M, \Omega_{k+1}(\E))
\end{equation}
and can thus be decomposed as a sum $\X = \X_- + \X_+$ onto each summand of \eqref{eq:XX}. The operator $\X_+$ is elliptic and has finite-dimensional kernel (when $M$ is compact) whose elements are called \emph{twisted conformal Killing tensors}. Moreover, the mapping property \eqref{eq:XX} ensures that $\X$ maps even (resp. odd) sections to odd (resp. even) sections. Elements in the kernel of $\X$ are \emph{flow-invariant}; equivalently, they have the property of invariance under parallel transport along geodesic flowlines.

\subsection{Examples. Link with Killing forms} 

\label{ssection:examples}

\subsubsection{Tautological section}

The tautological section, defined by $s(x,v) := v$, is a flow-invariant section of the pull-back bundle $\pi^*TM$ over $SM$ of Fourier degree 1. Flow-invariance is understood as above in the sense that $\X s = 0$, where $\X = (\pi^*\nabla)_X$ and $\nabla = \nabla^{TM}$ is the Levi-Civita connection on $TM$.
Equivalently, $s$ corresponds to the identity endomorphism $\Id_{TM}$, viewed as a section of $\Sym^1(T^*M)\otimes TM$, via the mapping \eqref{identif}, namely, $\pi^*(\Id_{TM})_{(x,v)} := \Id_{T_xM}(v) = v$. We use the notation $\Sym^1(T^*M)$ to insist on the fact that one could consider more general objects $\phi$ in $\Sym^p(T^*M) \otimes TM$ as in \eqref{identif}, and then the mapping to the unit tangent bundle would yield a section $(x,v) \mapsto \phi_x(v,...,v) \in \pi^*TM$.
 

\subsubsection{Normal bundle}
The normal bundle $\mc{V}$ on $SM$ is naturally identified with a subbundle of $\pi^*(TM)$ of codimension 1: as already mentioned, it is in fact the subbundle $s^\perp$ orthogonal to the tautological section $s$. Any section $f \in C^\infty(SM,\mc{V})$ which has Fourier degree 1 as a section of $\pi^*(TM)$ corresponds to an endomorphism $\phi$ of $TM$ via the above identification $f_{(x,v)} = \pi^*(\phi)_{(x,v)} = \phi_x(v)$, which further satisfies $g(\phi(v),v)=0$ for every $v\in SM$, i.e. it is skew-symmetric. 

\subsubsection{Exterior forms}
More generally, if $\omega$ is a $(p+1)$-form on $M$, it can be viewed as a $p$-form on $SM$ taking values in the normal bundle $\mc{V}$ by defining $\pi^*\omega \in C^\infty(SM,\Lambda^{p}\mc{V})$ as $\pi^*\omega_{(x,v)} := v \lrcorner \omega_x$ (interior product with $v$). Conversely, a section of Fourier degree 1 of $\pi^*(\Lambda^{p}(T^*M))$ which takes values in the subbundle $\Lambda^{p}\mc{V}$ of  $\pi^*(\Lambda^{p}(T^*M))$, corresponds to a $(p+1)$-form on $M$. Indeed, if $\omega\in C^\infty(M,\Sym^1(T^*M)\otimes \Lambda^p(T^*M))$ has the property that $(\pi^*\omega)_{(x,v)}\in \Lambda^p(\mc{V})$ for every $v\in SM$, this just means that $\omega$ is totally skew-symmetric as a covariant $(p+1)$-tensor.

\subsubsection{Flow-invariance and Killing forms}
It can be easily checked that the flow-invari\-ance condition $\X \pi^*\omega = 0$ translates into $\nabla_v^{TM} \omega (v, \bullet, ..., \bullet) = 0$ for every $v \in TM$, that is, the covariant derivative of $\omega$ is totally skew-symmetric. Such a form is called a \emph{Killing form} on $M$. A typical example where such a situation occurs with $p=1$ is the case of a Kähler manifold $(M,g,J)$, where $J^2 = -\mathbbm{1}_{TM}$ is the almost-complex structure satisfying $\nabla^{TM} J = 0$. Since $J$ is skew-symmetric, it defines an element of $\Lambda^2(T^*M)$. Setting $f_{(x,v)} := \pi^*(J)_{(x,v)} = J_x v$, the above discussion shows that $f \in C^\infty(SM,\mc{V})$ is flow-invariant. In turn, by \S \ref{ssection:invariant}, this shows that the transitivity group $H$ \emph{cannot be} equal to $\mathrm{SO}(n-1)$ since it has to fix at least one invariant vector, that is, $H \leqslant \mathrm{SO}(n-2)$, hence showing that the frame flow is not ergodic.


%
%
%

\subsection{Bounding the degree via the Pestov identity}

\label{ssection:pestov}

We will now explain the last steps in the proof of Theorem \ref{theorem:main1}. Assuming that the frame flow is not ergodic, Proposition \ref{prop:Brin1} implies that the transitivity group $H$ is a strict subgroup of $\mathrm{SO}(n-1)$, so by Theorem \ref{fis} there exists a non-vanishing section $f \in C^\infty(SM,\E)$, where $\E = \Lambda^p \mc{V}\subset \pi^*(\Lambda^{p}(T^*M))$ (with $p=1,2,3$) or $\E = \mathrm{Sym}^2 \mc{V}\subset \pi^*( \mathrm{Sym}^2(T^*M))$, satisfying $\X f =0$. The ultimate goal is to prove that $f$ is of degree $1$ under some pinching condition. The key identity for that is the twisted Pestov identity, originally stated in \cite[Proposition 3.5]{Guillarmou-Paternain-Salo-Uhlmann-16}, and restated in its present form in \cite[Lemma 2.3]{Cekic-Lefeuvre-Moroianu-Semmelmann-21}:

\begin{lemma}[Twisted Pestov identity]
\label{lemma:pestov}
Let $(M,g)$ be an $n$-dimensional compact Riemannian manifold and ${E}$ a Euclidean vector bundle over $M$, endowed with a metric connection $\nabla^{E}$.
If $\E$ denotes the pull-back of ${E}$ to $SM$, the following identity holds for all $k \in \Z_{\geq 0}$, and $u \in C^\infty(M,\Omega_k(\E))$:
\begin{equation}
\label{equation:local-pestov}
\frac{(n+k-2)(n+2k-4)}{n+k-3} \|\X_-u\|^2_{L^2} -   \frac{k(n+2k)}{k+1} \|\X_+u\|^2_{L^2} + \|Z u\|^2_{L^2}  = Q_k^{E}(u,u),
\end{equation}
where $Z$ is a first order differential operator which we do not make explicit and $Q_k^{E}$ is a symmetric bilinear form explicitly defined in terms of the Riemannian curvature of $(M,g)$ and the curvature tensor of $\nabla^{E}$.
\end{lemma}

Using the Cauchy-Schwarz inequality, one can show that when $(M,g)$ has sectional curvature bounded from above by $-\delta < 0$,
\begin{equation}
\label{qk}
Q_k^{E}(u,u)\leq (-\delta k^2+kq({E}))\|u\|^2_{L^2},
\end{equation}
where $q({E})$ only depends on the curvature tensor of $\nabla^{E}$. In particular, there exists an integer $k_0$ such that $Q_k^{E}(u,u)\leq 0$ for every $k\ge k_0$. 

Now, if $f \in C^\infty(SM,\E)$ satisfies $\X f=0$, we write $f=\sum_{k\ge 0}f_k$ with $f_k\in  C^\infty(M,\Omega_k(\E))$ which satisfy  $\|f_k\|_{H^1}\to 0$ as $k\to\infty$. Moreover, by \eqref{eq:XX} writing $\X=\X_++\X_-$ gives
\begin{equation}\label{pm}
	\X_+ f_{k}+\X_-f_{k+2}=0,\qquad \mbox{for every }k\ge 0. 
\end{equation}

Applying \eqref{equation:local-pestov} to $u=f_{k}$ we get for every $k> k_0$:
$$\frac{(n+k-2)(n+2k-4)}{n+k-3} \|\X_-f_{k}\|^2_{L^2} \le   \frac{k(n+2k)}{k+1} \|\X_-f_{k+2}\|^2_{L^2}.$$
In particular, this implies $\|\X_-f_{k}\|^2_{L^2} \le  \|\X_-f_{k+2}\|^2_{L^2}$ for every $k\ge k_0$. On the other hand, $\|\X_-f_{k}\|^2_{L^2}$ tends to $0$ as $k\to\infty$ because $\|\X_-f_k\|^2_{L^2} \leq \|\X f_k\|^2_{L^2} \leq \|f_k\|^2_{H^1} \to 0$ by smoothness of $f$, so $\X_-f_{k}=0$ for $k\ge k_0$. By \eqref{pm} we then also have $\X_+f_{k}=0$ for $k\ge 2+k_0$, so using \eqref{equation:local-pestov} and \eqref{qk} for $u:=f_k$ shows that $f_k=0$ for $k\ge 2+k_0$.
This gives the following result, originally proved in \cite[Theorem 4.1]{Guillarmou-Paternain-Salo-Uhlmann-16}:
\begin{corollary}\label{c42}
Every $f \in C^\infty(SM,\E)$ with $\X f=0$ has finite degree. 
\end{corollary}

Corollary \ref{c42} shows that if the frame flow of $(M,g)$ is not ergodic, the non-vanishing flow-invariant section $f \in C^\infty(SM,\E)$ given by Theorem \ref{fis} has finite degree. The idea is that under a suitable pinching hypothesis, one can show that the section $f$ has degree 1. By \S\ref{ssection:examples}, it defines a Killing form on $M$, the existence of which is obstructed either by negative curvature, or by its algebraic properties. We now explain the remaining arguments leading to Theorem \ref{theorem:main1} in cases (1)--(3)\footnote{Case (4) is slightly more technical and we refer to \cite[Section 4]{Cekic-Lefeuvre-Moroianu-Semmelmann-21} for further details.} of Theorem \ref{fis}.  

\begin{proof}[End of the proof of Theorem \ref{theorem:main1}] The proof is divided into three steps. \\

{\bf Step 1.} We first show by topological arguments that the section $f$ given by Theorem \ref{fis} is odd. For instance, in case (1), if non-zero, the restriction of the even part of $f$ to a fiber of $SM$ defines a constant length vector field on $\mathbb{S}^{n-1}$ of even degree, thus a polynomial map $\xi:\mathbb{S}^{n-1}\to \mathbb{S}^{n-1}$ satisfying $\xi(v)=\xi(-v)$ for every $v\in\mathbb{S}^{n-1}$. In particular, the topological degree of $\xi$ is even. On the other hand $\xi$ is homotopic to the identity via $\left[0, \frac{\pi}{2}\right] \ni t \mapsto \xi_t(v):=\cos(t)\xi(v)+\sin(t)v$, so its topological degree is 1, which is a contradiction. \\

{\bf Step 2.} Under some curvature pinching assumption, we then show that the degree of $f$ must be strictly smaller than $3$, hence precisely equal to $1$ by the first step. This is the key point of the proof, and is based on subtle estimates in the curvature term $Q_k^{E}$ appearing in the right-hand side of the twisted Pestov identity \eqref{equation:local-pestov}. In order to keep the discussion simple, we will only give the main idea. Decomposing the flow-invariant section $f = f_k + f_{k-2} + ... + f_1$, where $f_k \neq 0$ and $k$ is odd, and setting $u := f_k \neq 0$, we see by \eqref{eq:XX} that the flow-invariance $\X f = 0$ implies $\X_+ u = 0$, that is, $u$ is a twisted conformal Killing tensor. Applying \eqref{equation:local-pestov}, we thus obtain:
\begin{equation}
\label{equation:cont}
0 \leq \frac{(n+k-2)(n+2k-4)}{n+k-3} \|\X_-u\|^2 + \|Zu\|^2  = Q_k^{E}(u,u) \leq F(k,\delta)\|u\|^2,
\end{equation}
where $F(k,\delta)$ is the maximum of the symmetric bilinear form $Q_k^E$ on the unit sphere of $\Omega_k(\E)$. As \eqref{qk} indicates, it can be shown that for fixed $\delta$, the sequence $k \mapsto F(k,\delta)$ decreases to $-\infty$. Hence, there is a $k(\delta)$ such that for all $k \geq k(\delta)$, $F(k,\delta) < 0$. In turn, this implies by \eqref{equation:cont} that $u \equiv 0$ which contradicts the assumption that $u \neq 0$. Now, it can be checked that the function $\delta \mapsto k(\delta)$ is a \emph{decreasing} function, and that there exists $\delta(n) < 1$ such that for $\delta > \delta(n)$ sufficiently close to $1$, $k(\delta) < 3$. Hence, we conclude that whenever $\delta > \delta(n)$, if the twisted conformal Killing tensor $u$ in \eqref{equation:cont} is of degree $\geq 3$, it must vanish identically. This is a contradiction and thus $u$ is of degree $1$ by the first step. \\

{\bf Step 3.} Once we have established that the Fourier degree of the flow-invariant section $f \in C^\infty(SM,\Lambda^p \mc{V})$ is 1, following \S\ref{ssection:examples}, we know that $f$ defines a Killing $(p + 1)$-form on $M$, with $p \in \left\{1,2,3\right\}$, that is a $(p + 1)$-form $\omega$ such that $\nabla^{TM}_v \omega (v, \bullet, ..., \bullet) = 0$. The first two cases are ruled out by \cite{Barberis-Moroianu-Semmelmann-20}, while we use some ad-hoc arguments in order to show that such special Killing 4-forms vanish identically (see \cite[Lemma 3.13]{Cekic-Lefeuvre-Moroianu-Semmelmann-21}).
\end{proof}

\section{Open questions}

\label{section:questions}

In order to improve Theorem \ref{theorem:main1}, it is clear that a deeper algebraic understanding of the curvature term $Q_k^{E}$ appearing in the right-hand side of the twisted Pestov identity is missing. This is probably the main challenge for future work and numerical experiments might help studying this term more accurately. We list below some other open questions related to frame flow ergodicity:

\begin{itemize}
\item \textbf{Exponential mixing for the frame flow:} Using representation theory, exponential mixing was proved  by Moore \cite{Moore-87} on hyperbolic manifolds, and reproved recently in dimension $n=3$ by Guillarmou-Küster \cite{Guillarmou-Kuster-21} using semiclassical analysis. Rapid mixing of isometric extensions of diffeomorphisms was studied by Dolgopyat \cite{Dolgopyat-02} and more recently, by Siddiqi \cite{Siddiqi-19} in the case of flows. However, due to low regularity issues of the stable/unstable foliation, going beyond constant curvature is still an open question.
\item \textbf{$\mathrm{E}_7$-structure on $\Ss^{133}$:} As pointed out in Remark \ref{remark}, it is not known yet whether this possible structure actually does exist on $\Ss^{133}$. The main obstacle seems to be the computation of $\pi_{132}(\mathrm{E}_7)$ and this group is still unknown to the best of our knowledge.
\item \textbf{Existence of polynomial structures over spheres:} Whenever the frame flow is not ergodic, the arguments above give a flow-invariant section $f \in C^\infty(SM,\E)$ with finite Fourier degree. In turn, restricting to the sphere over a point $x \in M$, this entails the existence of a \emph{polynomial structure on the sphere} (for instance, a vector field whose coordinates are all homogeneous polynomials). These structures were studied by Wood \cite{Wood-68} but their complete classification is far from being understood.
\end{itemize}

\bibliographystyle{alpha}

\bibliography{Biblio}

\begin{thebibliography}{GPSU16}

\bibitem[Ada62]{Adams-62}
John~Frank Adams.
\newblock Vector fields on spheres.
\newblock {\em Ann. of Math. (2)}, 75:603--632, 1962.

\bibitem[Ami86]{Amirov-86}
Arif Amirov.
\newblock Existence and uniqueness theorems for the solution of an inverse
  problem for the transfer equation.
\newblock {\em Sibirsk. Mat. Zh.}, 27(6):3--20, 1986.

\bibitem[Ano67]{Anosov-67}
Dmitri~Viktorovitch Anosov.
\newblock Geodesic flows on closed {R}iemannian manifolds of negative
  curvature.
\newblock {\em Trudy Mat. Inst. Steklov.}, 90:209, 1967.

\bibitem[Ber60]{Berger-60-1}
Marcel Berger.
\newblock Pincement riemannien et pincement holomorphe.
\newblock {\em Ann. Scuola Norm. Sup. Pisa Cl. Sci. (3)}, 14:151--159, 1960.

\bibitem[BG80]{Brin-Gromov-80}
Michael Brin and Mikhael Gromov.
\newblock On the ergodicity of frame flows.
\newblock {\em Invent. Math.}, 60(1):1--7, 1980.

\bibitem[BK84]{Brin-Karcher-83}
Michael Brin and Hermann Karcher.
\newblock Frame flows on manifolds with pinched negative curvature.
\newblock {\em Compositio Math.}, 52(3):275--297, 1984.

\bibitem[BMS20]{Barberis-Moroianu-Semmelmann-20}
Mar\'{\i}a~Laura Barberis, Andrei Moroianu, and Uwe Semmelmann.
\newblock Generalized vector cross products and {K}illing forms on negatively
  curved manifolds.
\newblock {\em Geom. Dedicata}, 205:113--127, 2020.

\bibitem[BP74]{Brin-Pesin-74}
Michael~I. Brin and Yakov~B. Pesin.
\newblock Partially hyperbolic dynamical systems.
\newblock {\em Izv. Akad. Nauk SSSR Ser. Mat.}, 38:170--212, 1974.

\bibitem[BP03]{Burns-Pollicott-03}
Keith Burns and Mark Pollicott.
\newblock Stable ergodicity and frame flows.
\newblock {\em Geom. Dedicata}, 98:189--210, 2003.

\bibitem[Bri75]{Brin-75-1}
Michael Brin.
\newblock The topology of group extensions of {$C$}-systems.
\newblock {\em Mat. Zametki}, 18(3):453--465, 1975.

\bibitem[Bri82]{Brin-82}
Michael Brin.
\newblock Ergodic theory of frame flows.
\newblock In {\em Ergodic theory and dynamical systems, {II} ({C}ollege {P}ark,
  {M}d., 1979/1980)}, volume~21 of {\em Progr. Math.}, pages 163--183.
  Birkh\"{a}user, Boston, Mass., 1982.

\bibitem[BW10]{Burns-Wilkinson-10}
Keith Burns and Amie Wilkinson.
\newblock On the ergodicity of partially hyperbolic systems.
\newblock {\em Ann. of Math. (2)}, 171(1):451--489, 2010.

\bibitem[{\v{C}}C06]{Cadek-Crabb-06}
Martin {\v{C}}adek and Michael Crabb.
\newblock {$G$}-structures on spheres.
\newblock {\em Proc. London Math. Soc. (3)}, 93(3):791--816, 2006.

\bibitem[CL]{Cekic-Lefeuvre-21-1}
Mihajlo {Ceki{\'c}} and Thibault {Lefeuvre}.
\newblock {The Holonomy Inverse Problem}.
\newblock {\em arXiv:2105.06376}.

\bibitem[CLMS]{Cekic-Lefeuvre-Moroianu-Semmelmann-21}
Mihajlo {Ceki{\'c}}, Thibault {Lefeuvre}, Andrei {Moroianu}, and Uwe
  {Semmelmann}.
\newblock {On the ergodicity of the frame flow on even-dimensional manifolds}.
\newblock {\em arXiv:2111.14811}.

\bibitem[CS98]{Croke-Sharafutdinov-98}
Christopher~B. Croke and Vladimir~A. Sharafutdinov.
\newblock Spectral rigidity of a compact negatively curved manifold.
\newblock {\em Topology}, 37(6):1265--1273, 1998.

\bibitem[Dol02]{Dolgopyat-02}
Dmitry Dolgopyat.
\newblock On mixing properties of compact group extensions of hyperbolic
  systems.
\newblock {\em Israel J. Math.}, 130:157--205, 2002.

\bibitem[GK21]{Guillarmou-Kuster-21}
Colin Guillarmou and Benjamin K\"{u}ster.
\newblock Spectral theory of the frame flow on hyperbolic 3-manifolds.
\newblock {\em Ann. Henri Poincar\'{e}}, 22(11):3565--3617, 2021.

\bibitem[GPSU16]{Guillarmou-Paternain-Salo-Uhlmann-16}
Colin Guillarmou, Gabriel~P. Paternain, Mikko Salo, and Gunther Uhlmann.
\newblock The {X}-ray transform for connections in negative curvature.
\newblock {\em Comm. Math. Phys.}, 343(1):83--127, 2016.

\bibitem[Hop36]{Hopf-36}
Eberhard Hopf.
\newblock Fuchsian groups and ergodic theory.
\newblock {\em Trans. Amer. Math. Soc.}, 39(2):299--314, 1936.

\bibitem[HP06]{Hasselblatt-Pesin-06}
Boris Hasselblatt and Yakov Pesin.
\newblock Partially hyperbolic dynamical systems.
\newblock In {\em Handbook of dynamical systems. {V}ol. 1{B}}, pages 1--55.
  Elsevier B. V., Amsterdam, 2006.

\bibitem[{Lef}]{Lefeuvre-21}
Thibault {Lefeuvre}.
\newblock {Isometric extensions of Anosov flows via microlocal analysis}.
\newblock {\em arXiv:2112.05979}.

\bibitem[Leo71]{Leonard-71}
Peter Leonard.
\newblock {$G$}-structures on spheres.
\newblock {\em Trans. Amer. Math. Soc.}, 157:311--327, 1971.

\bibitem[Liv04]{Liverani-04}
Carlangelo Liverani.
\newblock On contact {A}nosov flows.
\newblock {\em Ann. of Math. (2)}, 159(3):1275--1312, 2004.

\bibitem[Moo87]{Moore-87}
Calvin~C. Moore.
\newblock Exponential decay of correlation coefficients for geodesic flows.
\newblock In {\em Group representations, ergodic theory, operator algebras, and
  mathematical physics ({B}erkeley, {C}alif., 1984)}, volume~6 of {\em Math.
  Sci. Res. Inst. Publ.}, pages 163--181. Springer, New York, 1987.

\bibitem[Muk75]{Mukhometov-75}
Ravil~Galatdinovich Mukhometov.
\newblock Inverse kinematic problem of seismic on the plane.
\newblock {\em Akad. Nauk. SSSR}, 6:243--252, 1975.

\bibitem[Muk81]{Mukhometov-81}
Ravil~Galatdinovich Mukhometov.
\newblock On a problem of reconstructing {R}iemannian metrics.
\newblock {\em Sibirsk. Mat. Zh.}, 22(3):119--135, 237, 1981.

\bibitem[Pes04]{Pesin-04}
Yakov~B. Pesin.
\newblock {\em Lectures on partial hyperbolicity and stable ergodicity}.
\newblock Zurich Lectures in Advanced Mathematics. European Mathematical
  Society (EMS), Z\"{u}rich, 2004.

\bibitem[PS88]{Pestov-Sharafutdinov-88}
Leonid~N. Pestov and Vladimir~A. Sharafutdinov.
\newblock Integral geometry of tensor fields on a manifold of negative
  curvature.
\newblock {\em Sibirsk. Mat. Zh.}, 29(3):114--130, 221, 1988.

\bibitem[PS00]{Pugh-Shub-00}
Charles Pugh and Michael Shub.
\newblock Stable ergodicity and julienne quasi-conformality.
\newblock {\em J. Eur. Math. Soc. (JEMS)}, 2(1):1--52, 2000.

\bibitem[PSU13]{Paternain-Salo-Uhlmann-13}
Gabriel~P. Paternain, Mikko Salo, and Gunther Uhlmann.
\newblock Tensor tomography on simple surfaces.
\newblock {\em Invent. Math.}, 193(1):229--247, 2013.

\bibitem[Sha94]{Sharafutdinov-94}
Vladimir~A. Sharafutdinov.
\newblock {\em Integral geometry of tensor fields}.
\newblock Inverse and Ill-posed Problems Series. VSP, Utrecht, 1994.

\bibitem[{Sid}]{Siddiqi-19}
Salman {Siddiqi}.
\newblock {Decay of correlations for certain isometric extensions of Anosov
  flows}.
\newblock {\em arXiv:1908.08550. {\rm To appear in} Ergodic Theory and
  Dynamical Systems}.

\bibitem[TZ]{Tsuji-Zhang-20}
Masato {Tsujii} and Zhiyuan {Zhang}.
\newblock {Smooth mixing Anosov flows in dimension three are exponential
  mixing}.
\newblock {\em arXiv:2006.04293}.

\bibitem[Woo68]{Wood-68}
Reginald Wood.
\newblock Polynomial maps from spheres to spheres.
\newblock {\em Invent. Math.}, 5:163--168, 1968.

\end{thebibliography}

\end{document}